# ON RANDOM MEASURES ON SPACES OF TRAJECTORIES AND STRONG AND WEAK SOLUTIONS OF STOCHASTIC EQUATIONS

**A. A. Dorogovtsev**                                                                    

We investigate stationary random measures on spaces of sequences or functions. A new definition of a strong solution of a stochastic equation is proposed. We prove that the existence of a weak solution in the ordinary sense is equivalent to the existence of a strong measure-valued solution.

## 1. Introduction

Let $(\mathfrak{X}, \rho)$ be a separable complete metric space and let $C = C(\mathbb{R}, \mathfrak{X})$ be the space of continuous functions acting from $\mathbb{R}$ to $\mathfrak{X}$ with the metric

$$\tilde{\rho}(f, g) \; = \; \sum_{k=1}^{\infty} \frac{1}{2^k} \varphi\Big( \max_{[-k, k]} \rho\big( f(t), g(t) \big) \Big).$$

Here, $\varphi(r) = r(1 + r)^{-1}$ for $r \geq 0$. With this metric, $C$ is a complete separable metric space. One of the main objects of investigation in this paper are random probability measures on $C$. Let us formulate the corresponding definition. Let $\mathfrak{A}$ denote the $\sigma$-algebra of Borel subsets of the space $C$ and let $(\Omega, \mathcal{F}, P)$ be a probability space.

**Definition 1.** *A random measure on $C$ is defined as a function $\mu \colon \Omega \times \mathfrak{A} \to [0; 1]$ satisfying the following conditions:*

*(i) for any fixed $\omega \in \Omega$, $\mu(\omega, \cdot)$ is a probability measure on $\mathfrak{A}$;*

*(ii) for any fixed $\Delta \in \mathfrak{A}$, $\mu(\cdot, \Delta)$ is a random variable.*

The investigation of random measures on a space of functions is motivated by the fact that measures of this type can be used as a tool for the description of the evolution of a mass on the original phase space $\mathfrak{X}$ [1, 2]. The example below illustrates this statement.

**Example 1.** Let $\mathfrak{X}$ be the space $\mathbb{R}^2$ with Euclidean metric. Consider the dynamical systems in $\mathfrak{X}$ defined by the differential equations

$$dx_1(t) = 0 dt, \quad dx_2(t) = T x_2(t) dt,$$

Institute of Mathematics, Ukrainian Academy of Sciences, Kiev.

where $T$ is a matrix of the form

$$T = \begin{pmatrix} 0 & -1 \\ 1 & 0 \end{pmatrix}.$$

The motion of space particles described by the first equation is trivial, i.e., all particles remain fixed. The motion of particles described by the second equation is a rotation around the origin of coordinates with constant velocity. On $\mathfrak{X}$, consider the initial mass distribution equal to the normal distribution with mean value zero and unit covariance matrix. Then the mass distribution remains unchanged both in the case of the motion of space particles according to the first equation and in the case of the motion of space particles according to the second equation. It is clear that the full picture of mass transfer is now given by finite-dimensional distributions, i.e., joint mass distributions for different finite collections of times. It is natural to replace this family of finite-dimensional distributions by a single measure on the space $C$.

In the present paper, we consider stationary random measures on $C$. To formulate the corresponding definition, we need the notion of random measures coinciding in distribution.

**Definition 2.** *We say that two random measures* $\mu_1$ *and* $\mu_2$ *on the space* $C$ *coincide in distribution if, for an arbitrary collection of sets* $\Delta_1, \ldots, \Delta_n \in \mathfrak{A}$, *the random vectors* $\big(\mu_1(\Delta_1), \ldots, \mu_1(\Delta_n)\big)$ *and* $\big(\mu_2(\Delta_1), \ldots, \mu_2(\Delta_n)\big)$ *have the same distribution.*

**Remark 1.** Definition 2 can be reformulated in terms of integrals of bounded measurable functions on $C$ with respect to the measures $\mu_1$ and $\mu_2$. The equivalent definition is formulated as follows:

**Definition 2′.** *We say that two random measures* $\mu_1$ *and* $\mu_2$ *on the space* $C$ *coincide in distribution if, for an arbitrary bounded measurable function* $f \colon C \to \mathbb{R}$, *the integrals*

$$\int_C f(u)\mu_i(du), \quad i = 1, 2,$$

*have the same distribution.*

**Remark 2.** The fact that the integral of a measurable function with respect to a random measure is a random variable is verified in a standard way, and, hence, we use it without proof.

Now let the mapping

$$C \ni u \to u_t \in C,$$

$$u_t(s) = u(s+t), \quad s \in \mathbb{R},$$

be the translation mapping on $C$. This is a measurable bijection on $C$. Therefore, for an arbitrary random measure $\mu$, its image $\mu^t$ under translation in $C$ is also a random measure.

Let us formulate the main definition of this paper.

**Definition 3.** *A random measure* $\mu$ *on* $C$ *is called stationary if, for any* $t \in \mathbb{R}$, *the random measures* $\mu$ *and* $\mu^t$ *are identically distributed.*

Random measures on $\mathbb{R}^{\mathbb{Z}}$, translations of random measures on $\mathbb{R}^{\mathbb{Z}}$, and stationary random measures on $\mathbb{R}^{\mathbb{Z}}$ are defined by analogy.

One of the aims of this paper is to describe the structure of stationary measures and construct nontrivial examples of such measures. It turns out that random measures on spaces of functions arise in a natural way in the course of the construction of solutions of stochastic differential and difference equations. In the case where there is no strong solution (i.e., a solution that is a functional of the original flow of random perturbations), a random measure can be used as its substitute. The paper consists of two parts. In the first part, we study the structure of stationary measures and give the corresponding examples. The second part is devoted to the use of random measures instead of a strong solution of a stochastic equation.

## 2. Stationary Random Measures on $C$

The theorem below provides examples of stationary random measures.

**Theorem 1.** *Suppose that there exists a group* $\left( T_t; t \in \mathbb{R} \right)$ *of measurable mappings on the probability space and, on the product of the spaces* $\Omega$ *and* $C$ *(with the* $\sigma$-*algebra spanned by the product of* $\sigma$-*algebras in* $\Omega$ *and* $C$), *there is a probability measure* $\tilde{P}$ *satisfying the following conditions:*

*(i) the mapping*

$$\Omega \times C \ni (\omega, u) \;\mapsto\; \left( T_t(\omega), u_t \right)$$

*preserves the measure* $\tilde{P}$ *for every* $t \in \mathbb{R}$;

*(ii) the projection of the measure* $\tilde{P}$ *onto* $\Omega$ *coincides with* $P$;

*(iii) the* $\sigma$-*algebra* $\mathcal{F}$ *is countably generated.*

*Then the family* $\left\{ \mu_{\omega}, \omega \in \Omega \right\}$ *of conditional measures of the measure* $\tilde{P}$ *with respect to* $\Omega$ *is a stationary random measure on* $C$.

**Remark 3.** The group $\left\{ T_t; t \in \mathbb{R} \right\}$ is additive with respect to the index $t$, i.e., $T_{t+s} = T_t \circ T_s$ for all $t$ and $s \in \mathbb{R}$. Moreover, the $\sigma$-algebra of random events $\mathcal{F}$ is assumed to be countably generated, which is a natural assumption guaranteeing the existence of a regular conditional probability.

**Proof.** Consider a measurable set $\Gamma$ in $\Omega \times C$. For an arbitrary $t \in \mathbb{R}$, let $\Gamma^t$ denote a translation of $\Gamma$, i.e., a set of the form $\Gamma^t = \left\{ (\omega, u_t) : (\omega, u) \in \Gamma \right\}$. Furthermore, for every $\omega \in \Omega$, we denote by $\Gamma_{\omega}$ a section of the set $\Gamma$, namely,

$$\Gamma_\omega = \{ u : (\omega, u) \in \Gamma \}.$$

By virtue of the definition of conditional measure, we have

$$\tilde{P}(\Gamma^t) = \int_\Omega \mu_\omega^{-t}(\Gamma_\omega) P(d\omega).$$

On the other hand,

$$\tilde{P}(\Gamma^t) = \tilde{P}\big(\{(\omega, u) : (T_{-t}(\omega), u_{-t}) \in \Gamma^t\}\big) = \tilde{P}\big(\{(\omega, u) : (T_{-t}(\omega), u) \in \Gamma\}\big)$$

$$= \int_\Omega \mu_\omega(\{ u : (T_{-t}(\omega), u) \in \Gamma \}) P(d\omega) = \int_\Omega \mu_{T_t(\omega)}(\{ u : (\omega, u) \in \Gamma \}) P \circ T_t(d\omega)$$

$$= \int_\Omega \mu_{T_t(\omega)}(\Gamma_\omega) P(d\omega).$$

In the last equality, we have used the fact that, according to the conditions of the theorem, the mappings $\{ T_t ; t \in \mathbb{R} \}$ preserve the measure $P$. Thus, for an arbitrary measurable set $\Gamma \subset \Omega \times C$, we have

$$\int_\Omega \mu_\omega^{-t}(\Gamma_\omega) P(d\omega) = \int_\Omega \mu_{T_t(\omega)}(\Gamma_\omega) P(d\omega).$$

Consequently,

$$\mu_\omega^{-t} = \mu_{T_t(\omega)} \pmod{P}.$$

However, the random measure $\{ \mu_{T_t(\omega)}, \omega \in \Omega \}$ has the same distribution as $\mu$ because $T_t$ preserves the measure $P$.

The theorem is proved.

Theorem 1 enables one to construct nontrivial examples of stationary random measures.

**Example 2.** Let $(\mathcal{Y}, d)$ be a complete metric separable space. In the product $\mathfrak{X} \times \mathcal{Y}$ with natural metric, we consider a two-component process $\{ (x_t, y_t) ; t \in \mathbb{R} \}$ stationary in the narrow sense and almost-surely continuous. Then the random measure $\mu$ defined by the relation

$$\mu(\Delta) = P \{ x \in \Delta / y \},$$

i.e., the conditional distribution of $x$ for known $y$, is a stationary random measure.

Theorem 1 admits an inversion in a certain sense. In other words, any random measure can be obtained as described above.

**Theorem 2.** *Let* $\mu$ *be a stationary random measure on* $C$. *Then there exists a new probability space* $(\Omega_1, \mathcal{F}_1, P_1)$, *a probability measure* $\tilde{P}_1$ *on* $\Omega \times C$, *and a family of mappings* $\{T_t\}$ *satisfying the conditions of Theorem 1 and such that the resulting stationary measure* $\mu_1$ *and* $\mu$ *are identically distributed.*

**Proof.** Consider the space $\mathfrak{M}$ of all probability measures on $C$ with the metric of weak convergence. Random measures described in Definition 1 are random elements in $\mathfrak{M}$. The translations $\{\mu^t; t \in \mathbb{R}\}$ of a stationary random measure $\mu$ form a stationary (in the narrow sense) random process in $\mathfrak{M}$. Let $(\Omega_1, \mathcal{F}_1, P_1)$ be the space $\mathfrak{M}^{\mathbb{R}}$ with the $\sigma$-algebra generated by cylindrical sets and the measure determined by the distribution of the process $\{\mu^t; t \in \mathbb{R}\}$. We define the probability measure $\tilde{P}_1$ on the product $\Omega_1 \times C$ as follows: For a measurable $\Gamma \subset \Omega_1 \times C$, we set

$$\tilde{P}_1(\Gamma) \ = \ \int_{\Omega_1} \omega_0(\Gamma_\omega) P_1(d\omega).$$

Here, as before, $\Gamma_\omega$ is the section of $\Gamma$ and $\omega_0$ is the value of the function $\omega : \mathbb{R} \to \mathfrak{M}$ at the point $0$. We define a family $\{T_t; t \in \mathbb{R}\}$ in the following way:

$$T_t(\omega, u) \ = \ \big(\omega(\,\cdot\, + t), u(\,\cdot\, + t)\big).$$

Let us verify that $T_t$ preserves the measure $\tilde{P}_1$. Indeed,

$$\tilde{P}_1 \circ T_t^{-1}(\Gamma) \ = \ \tilde{P}_1\big(\{(\omega, u) : T_t(\omega, u) \in \Gamma\}\big) \ = \ \tilde{P}_1\big(\{(\omega, u) : (\omega(\,\cdot\, + t), u(\,\cdot\, + t)) \in \Gamma\}\big)$$

$$= \ \int_{\Omega_1} \omega_0(\{u : (\omega(\,\cdot\, + t), u(\,\cdot\, + t)) \in \Gamma\}) P_1(d\omega) \ = \ \int_{\Omega_1} \omega_{-t}(\{u : (\omega, u(\,\cdot\, + t)) \in \Gamma\}) P_1^t(d\omega)$$

$$= \ \int_{\Omega_1} \omega_0(\{u : (\omega, u) \in \Gamma\}) P_1(d\omega) \ = \ \tilde{P}_1(\Gamma).$$

The fact that the random measure $\omega \in \Omega_1$ has the same distribution as $\mu$ now follows from the definitions of $\tilde{P}_1$ and $P_1$.

The theorem is proved.

For what follows, we need the definition of a random measure on $C$ consistent with a given filtration. Let $\{\mathcal{F}_t; t \geq 0\}$ be an increasing flow of $\sigma$-algebras on the original probability space $(\Omega, \mathcal{F}, P)$, $\mathcal{F}_t \subset \mathcal{F}$, $t \geq 0$. For any $t \geq 0$, let $\mathfrak{A}_t$ denote the $\sigma$-algebra of subsets of the space $C$ generated by coordinate functionals at points up to time $t$ inclusive.

**Definition 3.** *A random measure* $\mu$ *is called* $(\mathfrak{A}_t, \mathcal{F}_t)_{t \geq 0}$-*consistent if, for any* $t \geq 0$ *and* $\Delta \in \mathfrak{A}_t$, *the value of* $\mu(\Delta)$ *is an* $\mathcal{F}_t$-*measurable random variable.*

Note that, similarly to the case of stationary random measures, the above definition can be reformulated in terms of integrals of measurable bounded functions on $C$. A characterization of consistent random measures absolutely continuous with respect to a certain deterministic random measure is given in [3].

## 3. Strong and Weak Solutions of Stochastic Differential Equations

Consider the following Cauchy problem for stochastic differential equations in $\mathbb{R}$ on the interval $[0;1]$:

$$dx(t) = a(x_t, t)dt + b(x_t, t)dw(t),$$

$$x(0) = x_0. \tag{1}$$

Here, $\{w(t); t \in [0;1]\}$ is a Wiener process, $a$ and $b$ are measurable nonanticipating functionals on $C([0;1])$ [4], and, for every $t$, $x_t$ is the trajectory of a solution terminated at time $t$:

$$x_t(s) = x_{t \wedge s}, \quad s \in [0;1].$$

Recall the definitions of weak and small solutions of (1) (see, e.g., [4]).

**Definition 4.** *A random process* $\{x(t); t \in [0;1]\}$ *with continuous trajectories is called a strong solution of (1) if, for any* $t \in [0;1]$, *the random variable* $x(t)$ *is measurable with respect to the* $\sigma$-*algebra* $\sigma(w(s); s \le t)$, *and the integral relation corresponding to (1) is satisfied.*

**Definition 5.** *A pair of random processes with continuous trajectories* $\{x(t), w(t); t \in [0;1]\}$ *defined on a certain probability space* $(\Omega, \mathcal{F}, P)$ *with a selected flow of* $\sigma$-*algebras* $\{\mathcal{F}_t; t \in [0;1]\}$ *is called a weak solution of (1) if the following conditions are satisfied:*

(i) $w$ *is a* $\{\mathcal{F}_t; t \in [0;1]\}$-*Wiener process;*

(ii) *for any* $t \in [0;1]$, *the random variable* $x(t)$ *is measurable with respect to* $\mathcal{F}_t$, *and an integral analog of (1) is satisfied.*

It is known [5] that, in certain cases, the Cauchy problem (1) has a weak solution but does not have a strong solution. The advantage of a strong solution is that it is a functional of the original Wiener process. This property plays an important role, e.g., in the simulation of a solution or in the construction of approximations to a solution. Therefore, in the case where a strong solution does not exist, but there is a weak solution, it is desirable to have an object that is a functional of a given Wiener process and that can be used instead of a strong solution. The role of such an object may be played by a random measure on the space of continuous functions that is consistent with the flow of $\sigma$-algebras generated by the Wiener process. Consider one more situation where solutions measurable with respect to the original process describing random perturbations may not exist.

Let $\varphi : \mathbb{R}^2 \to \mathbb{R}$ be a measurable function and let $\{\xi_n; n \in \mathbb{Z}\}$ be a sequence of independent identically distributed random variables. Consider the recurrence equation

$$x_{n+1} = \varphi(x_n, \xi_{n+1}), \quad n \in \mathbb{Z}. \tag{2}$$

In certain cases, Eq. (2) has a stationary solution $\{x_n; n \in \mathbb{Z}\}$ stationarily related to $\{\xi_n; n \in \mathbb{Z}\}$, i.e., such that $\{(x_n, \xi_n); n \in \mathbb{Z}\}$ is a sequence stationary in the narrow sense. In [6], such solutions were obtained as the weak limits of solutions of the Cauchy problem for Eq. (2). In [7], a constructive method was used for finding solutions of this type for functions $\varphi$ close to linear ones (in the case of an infinite-dimensional phase space). The example below shows that a stationary solution of Eq. (2) is not always a strong solution.

**Example 3.** Let $\{\xi_n; n \in \mathbb{Z}\}$ be independent and uniformly distributed on $[0; 1]$ and let

$$\varphi(x, y) = \{x + y\}, \quad x, y \in \mathbb{R},$$

where $\{\cdot\}$ denotes the fractional part of a number. For $n \in \mathbb{Z}$, we set

$$\mathcal{F}_n = \sigma(\xi_k; k \leq n).$$

We now show that Eq. (2) does not have solutions $\{x_n; n \in \mathbb{Z}\}$ stationary in the narrow sense and such that, for every $n \in \mathbb{Z}$, $x_n$ is measurable with respect to $\mathcal{F}_n$. Assume the contrary. Let $\{x_n; n \in \mathbb{Z}\}$ be a solution that possesses the properties indicated. Consider the new random variables

$$y_n = e^{i2\pi x_n}, \quad n \in \mathbb{Z}.$$

Then

$$\forall n \in \mathbb{Z}: \quad y_{n+1} = y_n \cdot e^{i2\pi \xi_{n+1}}.$$

The sequence $\{y_n; n \in \mathbb{Z}\}$ is stationary and consistent with the sequence of $\sigma$-algebras $\{\mathcal{F}_n; n \in \mathbb{Z}\}$. Therefore,

$$\forall n \in \mathbb{Z}: \quad \mathrm{M} y_{n+1} = \mathrm{M} y_n \cdot \mathrm{M}\, e^{i2\pi \xi_{n+1}} = 0.$$

Furthermore, for fixed $n \in \mathbb{Z}$, we have

$$y_n = \lim_{m \to -\infty} \mathrm{M}(y_n / \sigma(\xi_m, \xi_{m+1}, \dots, \xi_n)) \pmod{P}.$$

However, it follows from the consistency of the sequence $\{y_n; n \in \mathbb{Z}\}$ that

$$\mathrm{M}(y_n / \sigma(\xi_m, \xi_{m+1}, \dots, \xi_n)) = \mathrm{M} y_{m-1} \cdot e^{i2\pi(\xi_m + \dots + \xi_n)} = 0.$$

Thus,

$$\forall n \in \mathbb{Z}: \quad y_n = 0 \quad (\text{mod } P),$$

which contradicts the equality

$$y_n = e^{i2\pi x_n}.$$

The contradiction obtained shows that the assumption of the existence of a stationary consistent solution (i.e., a strong solution) is not satisfied in the case considered. The reasoning presented above is a modified version of Tsirel'son's example [8]. Note that a stationary solution $\{x_n; n \in \mathbb{Z}\}$ stationarily related to $\{\xi_n; n \in \mathbb{Z}\}$ and such that, for every $n \in \mathbb{Z}$, $x_n$ and $\{\xi_k; k \geq n+1\}$ are independent exists and is unique. Indeed, let $\eta$ be a random variable uniformly distributed on $[0; 1]$ and independent of $\{\xi_n; n \in \mathbb{Z}\}$. We set $x_0 = \eta$ and define $x_n$ by the recurrence relation (2) for $n \geq 1$ and by the relation

$$x_{n-1} = \{x_n - \xi_n\}$$

for $n \leq -1$. Let us show that $\{x_n; n \in \mathbb{Z}\}$ is a stationary solution of Eq. (2) stationarily related to $\{\xi_n; n \in \mathbb{Z}\}$. For arbitrary $m \in \mathbb{Z}$ and $k \geq 0$, consider the distribution of $\{(x_m, \xi_m), (x_{m+1}, \xi_{m+1}), \dots, (x_{m+k}, \xi_{m+k})\}$. First, note that, for every $m \in \mathbb{Z}$, $x_m$ is uniformly distributed on $[0; 1]$. This can easily be verified by induction. Therefore, for any $\alpha \in [0; 1]$, we have

$$\mathrm{P}\{x_1 < \alpha/\xi\} = \mathrm{P}\{\{x_0 + \xi_1\} < \alpha/\xi\} = \mathrm{P}\{x_0 < \alpha\} = \mathrm{P}\{x_1 < \alpha\},$$

i.e., $x_1$ is independent of the sequence $\xi$. By analogy, we conclude that, for every $m \in \mathbb{Z}$, $x_m$ does not depend on $\xi$. Consequently, the distribution of $\{(x_m, \xi_m), \dots, (x_{m+k}, \xi_{m+k})\}$ coincides with the distribution of $\{(x_0, \xi_0), \dots, (x_k, \xi_k)\}$. Also note that, for any other solution $\{y_n; n \in \mathbb{Z}\}$ stationarily related to $\{\xi_n; n \in \mathbb{Z}\}$ and independent of the "future" $\xi$, the distributions of the sequences $(x_n, \xi_n)$, $n \in \mathbb{Z}$, and $(y_n, \xi_n)$, $n \in \mathbb{Z}$, coincide.

The example presented above shows that the recurrence equation (2) does not always have a stationary strong solution. Let us formulate the definition of a stationary random measure that is a solution of Eq. (2).

**Definition 6.** *A stationary random measure* $\mu$ *on* $\mathbb{R}^{\mathbb{Z}}$ *is called a strong solution of Eq. (2) if the following conditions are satisfied:*

(i) $\mu$ *is consistent with the sequence of* $\sigma$-*algebras* $\{\mathcal{F}_n = \sigma(\xi_k; k \leq n), n \in \mathbb{Z}\}$;

(ii) *for any* $m \geq 1$, $n \in \mathbb{Z}$, $\lambda_1, \dots, \lambda_m$, *and* $\rho \in \mathbb{R}$, *one has*

$$\int_{\mathbb{R}^{\mathbb{Z}}} \exp\left(i \sum_{k=1}^{m} \lambda_k u_{n+k} + i\rho u_{n+m+1}\right) \mu(du) = \int_{\mathbb{R}^{\mathbb{Z}}} \exp\left(i \sum_{k=1}^{m} \lambda_k u_{n+k}\right) \exp\left(i\rho \varphi(u_{n+m}, \xi_{n+m+1})\right) \mu(du). \tag{3}$$

**Remark 4.** Relation (3) is an analog of the Hopf equation for a space-time solution of the Navier–Stokes equation [9].

**Theorem 3.** *Equation (2) has a stationary solution $\{x_n; n \in \mathbb{Z}\}$ stationarily related to $\xi$ and such that, for any $n \in \mathbb{Z}$, $x_n$ is independent of $\{\xi_k; k \geq n + 1\}$, if and only if there exists a stationary random measure that is a strong solution of Eq. (2). If a solution $x$ possessing the properties indicated is unique, then a stationary measure that is a strong solution is also unique.*

**Proof.** Let $\{x_n; n \in \mathbb{Z}\}$ be a solution of Eq. (2) possessing the properties described in the theorem. Consider the random measure $\mu$ on $\mathbb{R}^{\mathbb{Z}}$ defined by the relation

$$\mu(u: (u_n, u_{n+1}, \ldots, u_{n+m}) \in \Delta) = P\{(x_n, x_{n+1}, \ldots, x_{n+m}) \in \Delta / \xi\}.$$

Here, $n \in \mathbb{Z}$, $m \geq 1$, and $\Delta$ is a Borel subset of $\mathbb{R}^{m+1}$. Since the space $\mathbb{R}^{\mathbb{Z}}$ can be metrized so that it becomes complete and separable and the $\sigma$-algebra of its Borel subsets coincides with the $\sigma$-algebra generated by cylindrical sets, it follows from the standard theorems that there exists a regular version of the conditional probability that determines the measure $\mu$. The fact that $\mu$ is a stationary random measure is verified in an obvious way. Let us check the consistency of $\mu$. For $n \in \mathbb{Z}$, we set

$$\mathcal{G}_n = \sigma((x_k, \xi_k); k \leq n),$$

$$\mathcal{F}_n = \sigma(\xi_k; k \leq n).$$

By assumption, for every $n \in \mathbb{Z}$, the $\sigma$-algebras $\mathcal{G}_n$ and $\sigma(\xi_k; k > n)$ are independent. Therefore, for any integrable random variable $\zeta$ measurable with respect to $\xi$ and any $n \in \mathbb{Z}$, we have

$$M(\zeta / \mathcal{G}_n) = M(\zeta / \mathcal{F}_n).$$

To prove this relation, it suffices to verify it for polynomials of the random variables $\{\xi_k; k \in \mathbb{Z}\}$, assuming that $\xi_0$ has all moments, and then pass to the limit. For a fixed $n \in \mathbb{Z}$, we now consider a bounded cylindrical function $\varphi$ on $\mathbb{R}^{\mathbb{Z}}$ dependent on coordinates whose indices do not exceed $n$. Then, for the random variable

$$\alpha = \int_{\mathbb{R}^{\mathbb{Z}}} \varphi(u) \mu(du) = M(\varphi(x)/\xi),$$

the following relations are true:

$$M\alpha\zeta = M\varphi(x)\zeta = MM(\varphi(x)\zeta / \mathcal{G}_n) = M\varphi(x)M(\zeta / \mathcal{G}_n) = M\varphi(x)M(\zeta / \mathcal{F}_n).$$

Here, as above, $\zeta$ is an arbitrary bounded random variable measurable with respect to $\xi$. Consequently,

$$\alpha = M(\varphi(x)/\mathcal{F}_n),$$

i.e., $\alpha$ is measurable with respect to $\mathcal{F}_n$. Thus, the random measure $\mu$ is consistent, i.e., it is a strong solution of Eq. (2).

Let us verify that the existence of a stationary measure $\mu$ that is a strong solution yields the existence of a stationary solution $\{x_n; n \in \mathbb{Z}\}$ (possibly on an extension of the original probability space). Consider a random sequence $x$ such that, for any measurable $\Delta \subset \mathbb{R}^{\mathbb{Z}}$, one has

$$P\{x \in \Delta/\xi\} = \mu(\Delta).$$

Note that, for the construction of $x$, one may need to extend the original probability space. The condition imposed on $\mu$ in the theorem yields

$$\forall n \in \mathbb{Z}, \quad \rho \in \mathbb{R}: \quad \mu\{u: \exp i\rho u_{n+1} = \exp i\rho \varphi(u_n, \xi_{n+1})\} = 1.$$

Therefore,

$$\forall n \in \mathbb{Z}: \quad \mu\{u: u_{n+1} = \varphi(u_n, \xi_{n+1})\} = 1.$$

Consequently,

$$\forall n \in \mathbb{Z}: \quad P\{x_{n+1} = \varphi(x_n, \xi_{n+1})\} = M\mu\{u: u_{n+1} = \varphi(u_n, \xi_{n+1})\} = 1.$$

The independence of the "past" $x$ and the "future" $\xi$ of the sequence follows from the relations

$$\forall k \geq 0, \ n \in \mathbb{Z}: \quad P\{(x_{n-k}, \dots, x_n) \in \Delta\} = M\mu(u: (u_{n-k}, \dots, u_n) \in \Delta),$$

$$P\{(x_{n-k}, \dots, x_n) \in \Delta/\xi_{n+1}, \dots\} = M\{\mathbb{1}_{\{(x_{n-k}, \dots, x_n) \in \Delta\}}/\xi_{n+1}, \dots\}$$

$$= M\{M\{\mathbb{1}_{\{(x_{n-k}, \dots, x_n) \in \Delta\}}/\xi\}/\xi_{n+1}, \dots\}$$

$$= M\{\mu(u: (u_{n-k}, \dots, u_n) \in \Delta)/\xi_{n+1}, \dots\}$$

$$= M\mu(u: (u_{n-k}, \dots, u_n) \in \Delta),$$

where $\Delta$ is an arbitrary measurable subset of $\mathbb{R}^{k+1}$.

The uniqueness of a strong solution is verified by analogy with the proof of Theorem 2.

The theorem is proved.